\documentclass[11pt, a4paper]{amsart}
\usepackage{amsfonts}

\usepackage[dvips]{epsfig}
\usepackage{amsmath}
\usepackage{amsthm}
\usepackage{amsbsy}
\usepackage{amsgen}
\usepackage{amscd}
\usepackage{amsopn}
\usepackage{amstext}
\usepackage{amsxtra}

\newtheorem{theorem}{Theorem}[section]
\newtheorem{proposition}[theorem]{Proposition}
\newtheorem{lemma}[theorem]{Lemma}

\newtheorem{remark}[theorem]{Remark}

\input{tcilatex}

\begin{document}
\title{On the Excursion Sets of Spherical Gaussian Eigenfunctions}
\author{Domenico Marinucci and Igor Wigman}
\address{Department of Mathematics, University of Rome Tor Vergata}
\email{marinucc@mat.uniroma2.it}
\address{Institutionen f\"{o}r Matematik, Kungliga Tekniska H\"{o}gskolan
(KTH), Stockholm \\
currently at Cardiff University, Wales, UK}
\email{wigmani@cardiff.ac.uk}
\thanks{IW is supported by the Knut and Alice Wallenberg Foundation, grant
KAW.2005.0098}
\maketitle

\begin{abstract}
The high frequency behaviour for random eigenfunctions of the spherical
Laplacian has been recently the object of considerable interest, also
because of strong motivations arising from Physics and Cosmology. In this
paper, we are concerned with the high frequency behaviour of excursion sets;
in particular, we establish a Uniform Central Limit Theorem for the
empirical measure, i.e. the proportion of spherical surface where spherical
Gaussian eigenfunctions lie below a level $z$. Our proofs borrows some
techniques from the literature on stationary long memory processes; in
particular, we expand the empirical measure into Hermite polynomials, and
establish a uniform weak reduction principle, entailing that the asymptotic
behaviour is asymptotically dominated by a single term in the expansion. As
a result, we establish a functional central limit theorem; the limiting
process is fully degenerate.

\begin{itemize}
\item \textbf{Keywords and Phrases: }Gaussian Eigenfunctions, Excursion
Sets, Empirical Measure, High Energy Asymptotics, Functional Central Limit
Theorem

\item \textbf{AMS Classification: }60G60; 42C10, 60D05, 60B10

\item \textbf{PACS: }02.50Ey, 02.30Nw, 02.30Px
\end{itemize}
\end{abstract}

\section{Introduction}

\subsection{Background}

Let $\Delta _{S^{2}}$ be the usual Laplacian on the sphere $\mathbb{S}^{2}$,
i.e., given by
\begin{equation*}
\Delta _{S^{2}}=\frac{1}{\sin \theta }\frac{\partial }{\partial \theta }%
\left\{ \sin \theta \frac{\partial }{\partial \theta }\right\} +\frac{1}{%
\sin ^{2}\theta }\frac{\partial ^{2}}{\partial \varphi ^{2}}\text{ , }0\leq
\theta \leq \pi \text{ , }0\leq \varphi <2\pi
\end{equation*}%
in the usual spherical coordinates $\left( \theta ,\varphi \right) ,$ An
orthonormal family of eigenfunctions for the Laplacian is well-known to be
given by the (complex-valued) spherical harmonics $\left\{ Y_{lm}\right\}
_{m=-l,...,l},$ i.e.
\begin{equation*}
\Delta _{S^{2}}Y_{lm}=-l(l+1)Y_{lm}\text{ , }l=1,2,3,...,\int_{S^{2}}Y_{lm}(%
\theta ,\varphi )\overline{Y_{lm}}(\theta ,\varphi )\sin \theta d\theta
d\varphi =\delta _{m^{\prime }}^{m},
\end{equation*}%
where $\delta _{m^{\prime }}^{m}$ is the usual Kronecker delta. For
definiteness, we recall that in coordinates, the (complex-valued) spherical
harmonics are defined by%
\begin{eqnarray*}
Y_{lm}(\theta ,\varphi ) &=&\sqrt{\frac{2l+1}{4\pi }\frac{(l-m)!}{(l+m)!}}%
P_{lm}(\cos \theta )\exp (im\varphi )\text{,} \\
P_{lm}(x) &=&(-1)^{m}(1-x^{2})^{m/2}\frac{d^{l+m}}{dx^{l+m}}(x^{2}-1)^{l}%
\text{,}
\end{eqnarray*}%
$\left\{ P_{lm}(.)\right\} $ denoting associated Legendre polynomials. The $%
(2l+1)$ spherical harmonics $\left\{ Y_{lm}\right\} $ thus provide an
orthonormal system for the space $\mathcal{H}_{l}$ of eigenfunctions
corresponding to the eigenvalue $-l(l+1)$, and we have the expansion%
\begin{equation*}
L^{2}(\mathbb{S}^{2})=\bigoplus_{l=1}^{\infty }\mathcal{H}_{l}\text{.}
\end{equation*}%
We shall equip each of the spaces $\left\{ \mathcal{H}_{l}\right\} $ with a
rotation-invariant Gaussian probability measure, and focus on the asymptotic
behaviour of the random eigenfunctions $\left\{ f_{l}(.)\right\}
_{l=1,2,...},$ i.e.%
\begin{equation*}
f_{l}(x)=\sum_{m=-l}^{l}a_{lm}Y_{lm}(x)\text{ , }x\in \mathbb{S}^{2}\text{,}
\end{equation*}%
where the coefficients $\left\{ a_{lm}\right\} $ are mean zero (complex)
Gaussian random variables such that $Ea_{lm}=0$ , $Ea_{lm}\overline{a}%
_{lm^{\prime }}=(2l+1)^{-1}\delta _{m}^{m^{\prime }},$ $a_{lm}=(-1)^{m}%
\overline{a}_{l,-m},$ whence
\begin{equation*}
Ef_{l}(x)\equiv 0\text{ , }Ef_{l}^{2}(x)=(2l+1)E\left| a_{lm}\right| ^{2}=1%
\text{.}
\end{equation*}%
Equivalently, $f_{l}$ could be expanded in a real-valued basis; throughout
this paper, however, we stick to the complex-valued representation, which
has simpler symmetry properties. The random functions $\left\{
f_{l}(.)\right\} $ are easily seen to be isotropic, i.e. for all $g\in SO(3)$
(the rotation group on $\mathbb{R}^{3})$ we have
\begin{equation*}
f_{l}(gx)\overset{d}{=}f_{l}(x)\text{,}
\end{equation*}%
$\overset{d}{=}$ denoting identity in law, in the sense of stochastic
processes.

As well-known, the spherical Gaussian eigenfunctions emerge very naturally
also in the analysis of random fields, since they provide the Fourier
components in spectral representation expansions. The geometry of random
fields has been of course extensively studied in the mathematical
literature: we refer to the well-known book \cite{adlertaylor} for a
complete list of references. Several statistics for Gaussian fields and
their excursion sets have been investigated and many applications have also
been implemented on experimental data, especially in a cosmological\
framework, see for instance \cite{hikage} and the references therein.
Geometric features of random spherical harmonics have been studied for
instance by (\cite{Wig1},\cite{Wig2}) where sharp estimates for the variance
of the nodal length (i.e. the boundary of excursion sets for $z=0$) are
provided. In this paper, we focus on a different characterization of the
excursion sets, i.e. the behaviour of their empirical measure, to be defined
below.

\subsection{Main results}

Let us define the \emph{spherical harmonics} \emph{empirical measure} as
follows: for all $z\in (-\infty ,\infty ),$
\begin{equation}
\Phi _{l}(z):=\int_{S^{2}}\mathbb{I}(f_{l}(x)\leq z)dx\text{,}
\label{empmea}
\end{equation}%
where $\mathbb{I}(\cdot )$ is, as usual, the indicator function which takes
value one if the condition in the argument is satisfied, zero otherwise. In
words, the function $\Phi _{l}(z)$ provides the (random) measure of the set
where the eigenfunction lie below the value $z.$ For example, the value of $%
\Phi _{l}(z)$ at $z=0$ is related to the so-called Defect
\begin{equation*}
\mathcal{D}_{l}:=\mathrm{meas}\left( f_{l}^{-1}(0,\infty )\right) -\mathrm{%
meas}\left( f_{l}^{-1}(-\infty ,0)\right)
\end{equation*}%
by the straightforward transformation
\begin{equation*}
\mathcal{D}_{l}=4\pi -2\Phi _{l}(0).
\end{equation*}%
Of course, $4\pi -\Phi _{l}(z)$ provides the area of the excursion set $%
\mathcal{A}_{l}(z):=\left\{ x:f_{l}(x)>z\right\} .$ The empirical measure
can be viewed as a continuous analogue of the empirical distribution
function for sequences of random variables, and it is simply related to the
first Minkowski functional from convex geometry, see \cite{adlertaylor}.
Clearly, for all $z\in \mathbb{R}$,
\begin{equation*}
\mathbb{E}\left[ \Phi _{l}(z)\right] =4\pi \Phi (z),
\end{equation*}%
where $\Phi (.)$ is the cumulative distribution function of the standard
Gaussian. Our aim in this paper is to study the distribution of $\left\{
\Phi _{l}(z)-\Phi (z)\right\} ,$ uniformly w.r.t. $z$, asymptotically as $%
l\rightarrow \infty $. The following lemma deals with the variance of $\Phi
_{l}(z)$ as $l\rightarrow \infty $.

\begin{lemma}
\label{lem:var(Phi_l)} For every $z\in \mathbb{R}$,
\begin{equation*}
\func{Var}(\Phi _{l}(z)) = z^2 \phi(z)^2 \cdot \frac{1}{l} + O_{z}\left(
\frac{\log{l}}{l^2} \right),
\end{equation*}
where $\phi$ is the standard Gaussian probability density function.
\end{lemma}

In particular, for $z\neq 0$, Lemma \ref{lem:var(Phi_l)} gives the
asymptotic form of the variance as $l\rightarrow \infty $. In contrast, for $%
z=0$ (this case corresponds to the Defect), this yields only a $``o"$-bound
and one needs to work harder to obtain a precise estimate (see ~\cite{defect}%
). Lemma \ref{lem:var(Phi_l)} follows rather easily from the Hermite
expansion approach in section \ref{sec:Hermite exp}, and \eqref{eq:Var hl2},
and we omit a formal proof.

From Lemma \ref{lem:var(Phi_l)} it is then natural to normalize $\Phi _{l}(z)
$ and define the \emph{spherical harmonics empirical process} by
\begin{equation}
G_{l}(z):=\sqrt{l}\left[ \int_{S^{2}}\mathbb{I}\left( f_{l}(x)\leq z\right)
dx-\left\{ 4\pi \times \Phi (z)\right\} \right]   \label{emppro}
\end{equation}%
for $l=1,2,...$, $z\in (-\infty ,\infty )$. In a sense to be made rigorous
later, we are then led to investigate the high-frequency ergodicity
properties of Gaussian eigenfunctions (compare \cite{MaPeJMP}). Indeed, as $l
$ diverges, (\ref{emppro}) amounts to the distance between an average over a
trajectory $\left( \int_{S^{2}}\mathbb{I}(f_{l}(x)\leq z)dx\right) $ and an
average over the ensemble determined by the probability law
\begin{equation*}
\Phi (z)=\int_{\Omega }\mathbb{I}(f_{l}(x)\leq z)dP(\omega ).
\end{equation*}%
For $G_{l}(.),$ we establish a functional central limit theorem on a
suitable Skorohod space.

Instrumental for the main result is a uniform weak reduction principle,
where we show that the behaviour of the empirical process is uniformly
approximated by a quadratic polynomial in $\left\{ f_{l}(.)\right\} $ times
a simple deterministic function. As a consequence, the weak limit is a fully
degenerate random process, which can be realized as a single random variable
times the derivative of the Gaussian density. In the sequel, we use $%
\Rightarrow $ to denote weak convergence in the Skorohod space $D([-\infty
,\infty ]),$ equipped with the sup-norm, and $G_{\infty }(z)$ to label the
mean zero, degenerate Gaussian process
\begin{equation*}
G_{\infty }(z)=z\mathbb{\phi }(z)Z\text{ , }Z\sim N(0,1)\text{.}
\end{equation*}

\begin{theorem}
\label{thm:Gl->Ginf} (The Uniform Central Limit Theorem) As $l\rightarrow
\infty ,$%
\begin{equation*}
G_{l}(z)\Rightarrow G_{\infty }(z).
\end{equation*}
\end{theorem}

Considering the limiting expression for $G_{\infty }(z),$ we remark the
existence of a singularity for $z=0$; indeed, here $J_{q}(0)=0$, whence it
is easily seen that
\begin{equation}
G_{l}(0)\rightarrow _{p}0\text{.}  \label{gldef}
\end{equation}%
From (\ref{gldef}), we immediately gather that the weak reduction argument
no longer holds for the Defect (see (\cite{defect}))$:$ indeed, the limiting
behaviour of the sequence $\left\{ D_{l}\right\} $ is considerably more
complicated and will be considered in our forthcoming paper ~\cite{defect}.

The material in this paper attempts to provide a characterization as
complete as possible on the high energy behaviour of the sample realizations
of Gaussian eigenfunctions. A further natural question is to investigate
whether the weak convergence result may enjoy some form of uniformity across
different frequencies $l.$ For the following result, $\Rightarrow $ denotes
weak convergence in the Skorohod space $D\left\{ [-\infty ,\infty ]\times
\lbrack 0,1]\right\} ,$ and $W_{\infty }(z;r)$ denotes the mean zero
Gaussian process with covariance function
\begin{equation*}
E\left[ W_{\infty }(z_{1};r_{1})W_{\infty }(z_{2};r_{2})\right] =\left\{
r_{1}\wedge r_{2}\right\} z_{1}z_{2}\mathbb{\phi }(z_{1})\mathbb{\phi }%
(z_{2})\text{.}
\end{equation*}

\begin{theorem}
\label{thm:Wl->Winf} (Partial Sum Processes) As $L\rightarrow \infty $%
\begin{equation*}
W_{L}(z;r)\Rightarrow W_{\infty }(z;r)\text{.}
\end{equation*}
\end{theorem}

\begin{remark}
It is possible to envisage some statistical applications of Theorem \ref%
{thm:Wl->Winf}, for instance in the analysis of isotropic spherical random
fields (such as those related to the analysis of the Cosmic Microwave
Background radiation (CMB), see \cite{dodelson} for details). In such
circumstances, the Gaussian eigenfunctions can be identified with the
Fourier components of the field at frequency $l,$ and it is straightforward
to exploit Theorems \ref{thm:Gl->Ginf} and \ref{thm:Wl->Winf} to construct
tests for Gaussianity and isotropy, of Kolmogorov-Smirnov ($S_{L})$ or
Cramer-von Mises ($K_{L})$ type, by means of the statistics%
\begin{equation*}
S_{L}:=\sup_{z}\sup_{r}\left| W_{L}(z;r)\right| \text{ , }%
K_{L}:=\int_{0}^{1}\int_{-\infty }^{\infty }\left| W_{L}(z;r)\right| ^{2}dzdr
\end{equation*}%
where weak convergence ensures that%
\begin{equation*}
S_{L}\rightarrow _{d}\sup_{z}\sup_{r}\left| W_{\infty }(z;r)\right| \text{ ,
}K_{L}\rightarrow _{d}\int_{0}^{1}\int_{-\infty }^{\infty }\left| W_{\infty
}(z;r)\right| ^{2}dzdr\text{,}
\end{equation*}%
whence threshold values at given significance level can be easily tabulated.
We refer for instance to \cite{m2006}, \cite{MarPTRF} for more discussion on
these issues.
\end{remark}

\subsection{Plan of the Paper}

In Section 2 we discuss polynomial transforms of spherical random
eigenfunctions and we establish their asymptotic behaviour under
Gaussianity. Particular care is devoted to the derivation of exact
asymptotic rates. In Section 3 we discuss the structure of our proofs and
the significance of main results. In Section 4 and 5 we provide the proofs
of Theorems (\ref{thm:Gl->Ginf}) and (\ref{thm:Wl->Winf}). Some background
material and technical lemmas are collected in an Appendix.

\subsection{Acknowledgements}

We are grateful to Giovanni Peccati, Ze\'{e}v Rudnick and Mikhail Sodin for
many stimulating discussions on these and related issues. This work was
initiated during the workshop ``Random Fields and Stochastic Geometry" held
in February 2009, Banff, Canada, and the authors wish to use this
opportunity to thank the organizers, Robert Adler and Jonathan Taylor, for
the hospitality, and an excellent opportunity to exchange ideas and
challenges. A substantial part of this research was done during the second
author's visit to University of Rome ``Tor Vergata", and he would like to
acknowledge the extremely friendly and stimulating environment in the
institution, and the generous financial support.

\section{Polynomial transforms and their asymptotic behaviour}

Throughout this paper, we shall make extensive use of well-known results on
Hermite polynomials and higher-order moments of Gaussian random variables.
To fix notation, we first recall the standard definition of Hermite
polynomials%
\begin{equation*}
H_{q}(x)=(-1)^{q}\phi ^{-1}(x)\left\{ \frac{d^{q}}{dx^{q}}\phi (x)\right\}
\text{,}
\end{equation*}%
the first few examples being provided by%
\begin{eqnarray*}
H_{1}(x) &=&x\text{ , }H_{2}(x)=x^{2}-1\text{ , }H_{3}(x)=x^{3}-3x\text{ , }%
H_{4}(x)=x^{4}-6x^{2}+3\text{ , } \\
H_{5}(x) &=&x^{5}-10x^{3}+15x\text{ , }H_{6}(x)=x^{6}-15x^{4}+45x^{2}+15%
\text{ , }...\text{ }
\end{eqnarray*}%
We recall also the basic formulae
\begin{equation*}
E\left[ H_{q}(f_{l}(x))\right] \equiv 0
\end{equation*}%
and
\begin{equation*}
E\left[ H_{q_{1}}(f_{l}(x))H_{q_{2}}(f_{l}(y))\right] =\delta
_{q_{1}}^{q_{2}}q_{1}!\left\{ Ef_{l}(x)f_{l}(y)\right\} ^{q_{1}}=\delta
_{q_{1}}^{q_{2}}q_{1}!P_{l}(\left\langle x,y\right\rangle )^{q_{1}},
\end{equation*}%
where%
\begin{equation*}
P_{l}(\left\langle x,y\right\rangle ):=\frac{4\pi }{2l+1}%
\sum_{m=-l}^{l}Y_{lm}(x)\overline{Y}_{lm}(y)=E\left\{
f_{l}(x)f_{l}(y)\right\}
\end{equation*}%
are the Legendre polynomials and $\left\langle x,y\right\rangle =\cos
(d(x,y))$ where $d(x,y)$ is the usual spherical distance on $\mathbb{S}^{2}.$
For our arguments in the sequel, we shall need a few basic facts on the
asymptotic behavior of the polynomial transforms
\begin{equation}
h_{l;q}:=\int_{S^{2}}H_{q}(f_{l}(x))dx\text{ , for }q=1,2,3,....
\label{eq:hlq def}
\end{equation}%
Note first that $Eh_{l;q}=0$ , for all $l,q;$ also
\begin{equation*}
h_{l;1}=\int_{S^{2}}f_{l}(x)dx\equiv 0\text{ , for all }l\geq 1\text{;}
\end{equation*}%
more generally, and $h_{l;q}\equiv 0$, when $l$ and $q$ are both odd, by
symmetry conditions. The behaviour of quadratic transforms is discussed in %
\cite{MaPeJMP}, where non-Gaussian behaviour is also covered and uniform
bounds on the speed of convergence to the limiting distribution are
discussed, by means of a detailed analysis of the cumulants for spherical
harmonics coefficients. We write here a different proof based on Legendre
polynomials to introduce the techniques we shall use extensively, later in
this paper.

\begin{proposition}
(Proposition 2, \cite{MaPeJMP}) As $l\rightarrow \infty $ we have%
\begin{equation*}
\sqrt{l}h_{l;2}\rightarrow _{d}N(0,1)\text{.}
\end{equation*}
\end{proposition}

\begin{proof}
Note first that%
\begin{eqnarray*}
Var(h_{l;2}) &=&Var\left\{ \int_{S^{2}}H_{2}(f_{l}(x))dx\right\} =E\left\{
\int_{S^{2}}H_{2}(f_{l}(x))dx\right\} ^{2} \\
&=&\int_{S^{2}\times S^{2}}E\left[ H_{2}(f_{l}(x))H_{2}(f_{l}(y))\right] dxdy
\\
&=&\int_{0}^{1}P_{l}^{2}(t)dt=\frac{2}{2l+1}\text{.}
\end{eqnarray*}%
In view of results \cite{NouPe1, Noupe2}, to investigate the asymptotic
behaviour of the total variation distance $d_{TV}(\sqrt{l}h_{l;2},Z)$
between $\sqrt{l}h_{l;2}$ and a standard Gaussian random variable $Z\sim
N(0,1),$ it is enough to show that the normalized fourth-order cumulants
converge to zero, i.e., that
\begin{equation*}
\frac{cum_{4}(h_{l;2})}{\left\{ Var(h_{l;2})\right\} ^{2}}=o_{l\rightarrow
\infty }(1).
\end{equation*}%
Using the well-known Diagram Formula (see e.g. \cite{Surg}), we have%
\begin{equation*}
cum_{4}(h_{l;2})=48\int_{S^{2}\times S^{2}\times S^{2}\times
S^{2}}P_{l}(\left\langle x,y\right\rangle )P_{l}(\left\langle
y,z\right\rangle )P_{l}(\left\langle z,w\right\rangle )P_{l}(\left\langle
w,x\right\rangle )dxdydzdw\text{.}
\end{equation*}%
Now recall the reproducing kernel formula
\begin{equation*}
\int_{S^{2}}P_{l}(\left\langle x,y\right\rangle )P_{l}(\left\langle
y,z\right\rangle )dy=\frac{4\pi }{2l+1}P_{l}(\left\langle x,y\right\rangle )%
\text{,}
\end{equation*}%
whence,

\begin{equation*}
cum_{4}(h_{l;2})=48\left\{ \frac{4\pi }{2l+1}\right\} ^{2}\int_{S^{2}\times
S^{2}}P_{l}^{2}(\left\langle w,x\right\rangle )dxdw=48\frac{\left( 4\pi
\right) ^{4}}{\left( 2l+1\right) ^{3}}\text{,}
\end{equation*}%
We then have
\begin{equation*}
d_{TV}(\sqrt{l}h_{l;2},Z)=O\left( \frac{cum_{4}(h_{l;2})}{\left\{
Var(h_{l;2})\right\} ^{2}}\right) =O(l^{-1}),
\end{equation*}%
by the result mentioned above.
\end{proof}

The next Proposition covers the case $q=3$ and follows easily from results
in \cite{m2006, MarPTRF}.

\begin{proposition}
As $l\rightarrow \infty ,$ we have%
\begin{equation*}
lh_{l;3}\rightarrow _{d}N\left( 0,\frac{\sqrt{3}}{\pi ^{2}}\right) .
\end{equation*}%
\bigskip
\end{proposition}

\begin{proof}
Note first that%
\begin{eqnarray*}
h_{l;3}
&=&\int_{S^{2}}f_{l}^{3}(x)dx=%
\sum_{m_{1}m_{2}m_{3}}a_{lm_{1}}a_{lm_{2}}a_{lm_{3}}%
\int_{S^{2}}Y_{lm_{1}}Y_{lm_{2}}Y_{lm_{3}}dx \\
&=&\sum_{m_{1}m_{2}m_{3}}a_{lm_{1}}a_{lm_{2}}a_{lm_{3}}\left(
\begin{array}{ccc}
l & l & l \\
m_{1} & m_{2} & m_{3}%
\end{array}%
\right) \left(
\begin{array}{ccc}
l & l & l \\
0 & 0 & 0%
\end{array}%
\right) \sqrt{\frac{(2l+1)^{3}}{4\pi }} \\
&=&I_{lll}\left(
\begin{array}{ccc}
l & l & l \\
0 & 0 & 0%
\end{array}%
\right) \frac{1}{\sqrt{4\pi }}\text{,}
\end{eqnarray*}%
where%
\begin{equation*}
I_{lll}=\sum_{m_{1}m_{2}m_{3}}\frac{a_{lm_{1}}a_{lm_{2}}a_{lm_{3}}}{\sqrt{%
(2l+1)^{-3}}}\left(
\begin{array}{ccc}
l & l & l \\
m_{1} & m_{2} & m_{3}%
\end{array}%
\right) \text{,}
\end{equation*}%
is the normalized bispectrum (see \cite{Hu,m2006,MarPTRF} for definitions
and further discussions); the symbol on the right-hand side denotes the
so-called Wigner coefficients, for which we refer to \cite{VMK} and the
Appendix below. It was shown in \cite{MarPTRF} that
\begin{equation*}
EI_{lll}=0\text{ , }EI_{lll}^{2}=6\text{ and }I_{lll}\rightarrow _{d}N(0,6)%
\text{ as }l\rightarrow \infty \text{ .}
\end{equation*}%
Hence, in view of Lemma \ref{cgbou},
\begin{equation*}
\lim_{l\rightarrow \infty }l^{2}Eh_{l;3}^{2}=\frac{3}{2\pi }%
\lim_{l\rightarrow \infty }l^{2}\left(
\begin{array}{ccc}
l & l & l \\
0 & 0 & 0%
\end{array}%
\right) ^{2}=\frac{3}{2\pi }\times \frac{2}{\pi \sqrt{3}}=\frac{\sqrt{3}}{%
\pi ^{2}},
\end{equation*}%
as claimed.
\end{proof}

For the higher order terms, we have the following:

\begin{lemma}
\label{quattro} For all $q\geq 3,$ as $l\rightarrow \infty $ we have
\begin{equation}
Eh_{l;q}^{2}=q!\int_{S^{2}\times S^{2}}P_{l}(\left\langle x,y\right\rangle )
^{q} dxdy=O\left(\frac{\log l}{l^{2}}\right).  \label{rid}
\end{equation}
\end{lemma}

The proof of Lemma \ref{quattro} is given in Appendix \ref{apx:proof lem Leg
mom} for completeness. In fact, it is possible to obtain the asymptotics of
the LHS of \eqref{rid}:
\begin{equation*}
\int_{S^{2}\times S^{2}}P_{l}(\left\langle x,y\right\rangle ) ^{4} dxdy \sim
c_{4} \cdot \frac{\log {l}}{l^2}.
\end{equation*}
and for $q=3$ or $q\ge 5$ we have
\begin{equation*}
\int_{S^{2}\times S^{2}}P_{l}(\left\langle x,y\right\rangle ) ^{q} dxdy \sim
c_{q} \cdot \frac{1}{l^2},
\end{equation*}
for some constants $c_{q}$ (in particular, we get rid of the logarithm on
the RHS of \eqref{rid}). We use the precise estimates in our forthcoming
paper \cite{defect} in order to compute the variance of $\left\{
D_{l}\right\}$.

\section{On the proofs of the main results}

\label{sec:Hermite exp}

Our aim in this Section is to discuss heuristically the expansion of the
empirical measure and process (\ref{empmea})-(\ref{emppro}) into Hermite
polynomials. Recall first that $\left\{ f_{l}(x)\right\} $ is defined
pointwise as a Gaussian variable, whence the function $\mathbb{I}%
(f_{l}(x)\leq z)$ belongs for each $x,z$ to the $L^{2}$ space of square
integrable functions of Gaussian variables. In particular, we have that
\begin{equation*}
\mathbb{I}(f_{l}(x)\leq z)=\sum_{q=0}^{\infty }\frac{J_{q}(z)}{q!}%
H_{q}(f_{l}(x))\text{,}
\end{equation*}%
where the right hand side converges in the $L^{2}$ sense (i.e.
\begin{equation*}
\lim_{Q\rightarrow \infty }E\left\{ \sum_{q=Q}^{\infty }\frac{J_{q}(z)}{q!}%
H_{q}(f_{l}(x))\right\} ^{2}=0,
\end{equation*}%
uniformly w.r.t. $z,x$)$.$ It is possible to provide analytic expressions
for the coefficients $\left\{ J_{q}(.)\right\} ,$ indeed for $q\geq 1$%
\begin{eqnarray*}
J_{q}(z) &=&\int_{R}\mathbb{I}(u\leq z)H_{q}(u)\phi (u)du=\int_{-\infty
}^{z}(-1)^{q}\left\{ \frac{d^{q}}{du^{q}}\mathbb{\phi }(u)\right\} \phi
^{-1}(u)\phi (u)du \\
&=&(-1)^{q}\int_{-\infty }^{z}\left\{ \frac{d^{q}}{du^{q}}\mathbb{\phi }%
(u)\right\} du=(-1)^{q}\Phi ^{(q)}(z)=(-1)^{q}\mathbb{\phi }^{(q-1)}(z)\text{%
,}
\end{eqnarray*}%
where $\mathbb{\phi }(z):=d\Phi (z)/dz$ is the standard Gaussian density.
Hence for instance $J_{0}(z)=\Phi (z)$, $J_{1}(z)=-\mathbb{\phi }(z)$, $%
J_{2}(z)=-z\mathbb{\phi }(z)$, $J_{3}(z)=(1-z^{2})\mathbb{\phi }(z)$ and in
general
\begin{equation}  \label{eq:Jq(z)=-Hq(z) phi(z)}
J_{q}(z)=-H_{q}(z)\mathbb{\phi }(z).
\end{equation}
For $z=0$ it follows that $J_{1}(0)=-(2\pi )^{-1/2},$%
\begin{equation*}
J_{q}(0)=%
\begin{cases}
\frac{(-1)^{q}}{\sqrt{2\pi }}(q-2)!!\;\; & q=2k+1,\,k=1,2,... \\
0 & q=2k,\,k=1,2,.....%
\end{cases}%
\end{equation*}

We have, formally%
\begin{eqnarray*}
\int_{S^{2}}\mathbb{I}\left( f_{l}(x)\leq z\right) dx
&=&\int_{S^{2}}\sum_{q=0}^{\infty }\frac{J_{q}(z)}{q!}H_{q}(f_{l}(x))dx \\
&=&4\pi \times \Phi (z)+\int_{S^{2}}\sum_{q=2}^{\infty }\frac{J_{q}(z)}{q!}%
H_{q}(f_{l}(x))dx\text{,}
\end{eqnarray*}%
whence%
\begin{equation*}
G_{l}(z)=\sqrt{l}\int_{S^{2}}\sum_{q=2}^{\infty }\frac{J_{q}(z)}{q!}%
H_{q}(f_{l}(x))dx\text{,}
\end{equation*}%
and heuristically%
\begin{eqnarray}
G_{l}(z) &=&-\frac{z\mathbb{\phi }(z)}{2}\sqrt{l}%
\int_{S^{2}}H_{2}(f_{l}(x))dx  \notag \\
&&+\frac{(1-z^{2})\mathbb{\phi }(z)}{6}\sqrt{l}%
\int_{S^{2}}H_{3}(f_{l}(x))dx+... \\
&=&-\frac{z\mathbb{\phi }(z)}{2}\sqrt{l}h_{l;2}+\frac{(1-z^{2})\mathbb{\phi }%
(z)}{6}\sqrt{l}h_{l;3}-\frac{(3z-z^{3})\mathbb{\phi }(z)}{24}\sqrt{l}%
h_{l;4}+...,  \label{expa}
\end{eqnarray}%
where $h_{l;q}$ are as earlier in \eqref{eq:hlq def}. The meaning of the
exchange between the integral and the series is discussed below. From the
results in the previous section, we know that the variance of $\sqrt{l}%
h_{l;2}$ is given by
\begin{equation}
\func{Var}\left\{ \sqrt{l}h_{l;2}\right\} =\frac{2l}{2l+1},
\label{eq:Var hl2}
\end{equation}%
and the variances of the other terms are bounded by
\begin{equation*}
\func{Var}\left\{ \sqrt{l}h_{l;3}\right\} =O\left( \frac{1}{l}\right) \text{,%
}\;\func{Var}\left\{ \sqrt{l}h_{l;q}\right\} =O\left( \frac{\log l}{l}%
\right) \text{ , }q\geq 4.
\end{equation*}%
In fact, one can compute the precise asymptotic shape of the latter (see the
discussion after Lemma \ref{quattro}).

Our idea in the proofs to follow is hence to show that the term $\left\{
\sqrt{l}h_{l;2}\right\} $ is asymptotically dominant in the expansion of $%
G_{l}$, all the other summands being uniformly of smaller order. Note that,
as remarked in the Introduction, the first term in the expansion (\ref{expa}%
), ($q=2)$ vanishes for $z=0$; as a consequence, in the Defect statistic $%
D_{l}$ the quadratic term $2^{-1}z\mathbb{\phi }(z)\sqrt{l}h_{l;2} $
disappears. From the arguments below it is easy to realize that $D_{l}$ has
hence a lower order asymptotic variance than the empirical process; no
single term is asymptotically dominant and the determination of the limiting
behaviour becomes much more challenging. The asymptotic behaviour of the
sequence $\left\{ D_{l}\right\} $ is hence delayed to the forthcoming paper %
\cite{defect}.

\section{Proof of Theorem \ref{thm:Gl->Ginf}}

\begin{proof}
The key step in the proof is the following \emph{uniform weak reduction
principle: } there exist a constant $C$ such that for all $0<\varepsilon
\leq 1,$ as $l\rightarrow \infty ,$
\begin{equation*}
\Pr \left\{ \sup_{z}\left| G_{l}(z)-J_{2}(z)\sqrt{l}h_{l;2}dx\right|
>\varepsilon \right\} \leq \frac{C\log ^{6}l}{\varepsilon ^{3}l},
\end{equation*}%
which is the statement of Lemma \ref{chain}. We deduce from this that
\begin{equation*}
\sup_{z}\left| G_{l}(z)-J_{2}(z)\sqrt{l}h_{l;2}\right| =o_{p}(1),
\end{equation*}%
and the result follows immediately from
\begin{equation*}
\sqrt{l}h_{l;2}\rightarrow _{d}N(0,1)
\end{equation*}%
(see \cite{MaPeBer, MaPeJMP}).
\end{proof}

It then remains to prove Lemma \ref{chain}. Our approach will follow closely
the ideas by \cite{DelTaq}, which were developed in a rather different
context (the empirical process for long range dependent sequences on $%
\mathbb{R}).$ We will need the following notation:
\begin{equation}
S_{l}(z):=G_{l}(z)-J_{2}(z)\sqrt{l}h_{l;2}\text{,}  \label{important}
\end{equation}

\begin{equation*}
S_{l}(z_{1},z_{2}):=S_{l}(z_{2})-S_{l}(z_{1})\text{ },\text{ }z_{1}\leq
z_{2},
\end{equation*}%
and%
\begin{equation*}
J_{q}(z_{2},z_{1}):=J_{q}(z_{2})-J_{q}(z_{1})\text{ , }z_{1}\leq z_{2}\text{.%
}
\end{equation*}

\begin{lemma}
There exist a constant $C>0$ such that%
\begin{equation}
E\left[ S_{l}(z_{1},z_{2})\right] ^{2}\leq \frac{C\log l}{l}\left[ \Phi
(z_{2})-\Phi (z_{1})\right] \text{.}.  \label{sab1}
\end{equation}
\end{lemma}

\begin{proof}
For notational simplicity and without loss of generality we take $%
z_{1}=-\infty ,$ i.e. $S_{l}(z_{1},z_{2})=S_{l}(z_{2})=S_{l}(z);$ the
argument in the general case is identical. We have trivially
\begin{equation}  \label{eq:ES}
E\left\{ \int_{S^{2}}\left[ \mathbb{I}(f_{l}(x)\leq z)-\Phi (z)+\frac{z%
\mathbb{\phi }(z)}{2}H_{2}(f_{l}(x)) \right] dx\right\} ^{2}
\end{equation}
\begin{eqnarray*}
&=&E\left\{ \int_{S^{2}\times S^{2}}\left[ \mathbb{I}(f_{l}(x)\leq z)-\Phi
(z)+\frac{z\mathbb{\phi }(z)}{2}H_{2}(f_{l}(x)) \right] \right. \\
&&\left. \times \left[ \mathbb{I}(f_{l}(y)\leq z)-\Phi (z)+\frac{z\mathbb{%
\phi }(z)}{2}H_{2}(f_{l}(y)) \right] dxdy\right\}
\end{eqnarray*}
\begin{eqnarray*}
&=&\left\{ \int_{S^{2}\times S^{2}}E\left[ \mathbb{I}(f_{l}(x)\leq z)-\Phi
(z)+\frac{z\mathbb{\phi }(z)}{2}H_{2}(f_{l}(x)) \right] \right. \\
&&\left. \times \left[ \mathbb{I}(f_{l}(y)\leq z)-\Phi (z)+\frac{z\mathbb{%
\phi }(z)}{2} H_{2}(f_{l}(y)) \right] dxdy\right\} \text{.}
\end{eqnarray*}
Now%
\begin{equation}  \label{eq:E S as Leg}
\begin{split}
&E\left\{ \left[ \mathbb{I}(f_{l}(x)\leq z)-\Phi (z)+\frac{z\mathbb{\phi }
(z)}{2}H_{2}(f_{l}(x)) \right] \cdot \left[ \mathbb{I}(f_{l}(y)\leq z)-\Phi
(z)+\frac{z\mathbb{\phi }(z)}{2} H_{2}(f_{l}(y)) \right] \right\} \\
&=E\left[ \sum_{q_{1}=3}^{\infty }\frac{J_{q_{1}}(z)}{q_{1}!}%
H_{q_{1}}(f_{l}(x))\sum_{q_{2}=3}^{\infty }\frac{J_{q_{2}}(z)}{q_{2}!}%
H_{q_{2}}(f_{l}(y))\right] =\sum_{q=3}^{\infty }\frac{J_{q}^{2}(z)}{q!}%
P_{l}^{q}(\left\langle x,y\right\rangle )\text{.}
\end{split}%
\end{equation}
This result is standard, as $\sum_{q_{1}=3}^{\infty }\frac{J_{q_{1}}(z)}{%
q_{1}!}H_{q_{1}}(f_{l}(x))$ belongs to the Hilbert space of Gaussian
subordinated random variables $G(f_{l})$, with inner product $\left\langle
X,Y\right\rangle :=EXY.$ Plugging \eqref{eq:E S as Leg} into \eqref{eq:ES}
and using Lemma \ref{quattro}, we finally obtain
\begin{eqnarray*}
&&E\left\{ \int_{S^{2}}\left[ \mathbb{I}(f_{l}(x)\leq z)-\Phi (z)+\frac{z%
\mathbb{\phi }(z)}{2}H_{2}(f_{l}(x)) \right] dx\right\} ^{2} \\
&=&\int_{S^{2}\times S^{2}}\sum_{q=3}^{\infty }\frac{J_{q}^{2}(z)}{q!}%
P_{l}^{q}(\left\langle x,y\right\rangle )dxdy=\left\{ \sum_{q=3}^{\infty }%
\frac{J_{q}^{2}(z)}{q!}\right\} \times \left\{ \sup_{q\geq
3}\int_{S^{2}\times S^{2}}P_{l}^{q}(\left\langle x,y\right\rangle
)dxdy\right\} \\
&=&O\left(\frac{\log l}{l^{2}}\right)\text{.}
\end{eqnarray*}
since
\begin{equation*}
\sum_{q=3}^{\infty }\frac{J_{q}^{2}(z)}{q!}= \func{Var}\left(\mathbb{I}%
(f_{l}(x)\leq z)\right) = \Phi (z)(1-\Phi (z))<\infty \text{,}
\end{equation*}
is the variance of a Bernoulli random variable.
\end{proof}

\begin{lemma}
\label{chain} (Chaining argument) For all $0<\varepsilon \leq 1$%
\begin{equation}
\Pr \left\{ \sup_{z}\left| S_{l}(z)\right| >\varepsilon \right\} \leq \frac{C%
}{l}\left\{ \frac{\log ^{6}l}{\varepsilon ^{3}}\right\} \text{,}
\label{chain1}
\end{equation}%
and for $\varepsilon \geq 1$%
\begin{equation}
\Pr \left\{ \sup_{z}\left| S_{l}(z)\right| >\varepsilon \right\} \leq \frac{C%
}{l}\left\{ \varepsilon ^{-4/3}+\frac{\log ^{6}l}{\varepsilon ^{7/3}}%
\right\} \text{.}  \label{chain2}
\end{equation}
\end{lemma}

\begin{proof}
For (\ref{chain1}), we follow again closely \cite{DelTaq}. We define
\begin{equation*}
\Lambda (z)=\Phi (z)+\frac{1}{2}\int_{-\infty }^{z}\left| x^{2}-1\right|
\phi (x)dx
\end{equation*}%
\begin{equation*}
=\Phi (z)+\frac{1}{2\sqrt{2\pi }}\left\{ -z\exp (-\frac{z^{2}}{2})\mathbb{I}%
(z\leq -1)\right\}
\end{equation*}%
\begin{equation*}
+\frac{1}{2\sqrt{2\pi }}\left\{ \left[ z\exp (-\frac{z^{2}}{2})+\exp (-\frac{%
1}{2})\right] \mathbb{I}(-1<z\leq 1)+\left[ -z\exp (-\frac{z^{2}}{2})+4\exp
(-\frac{1}{2})\right] \mathbb{I}(z>1)\right\} \text{.}
\end{equation*}%
Clearly $\Lambda (z)$ is continuous, $\lim_{z\rightarrow -\infty }\Lambda
(z)=0,$ $\lim_{z\rightarrow \infty }\Lambda (z)=1.483943...=:\Lambda (\infty
).$ Also, for $\Phi (z_{1},z_{2})=\Phi (z_{2})-\Phi
(z_{1})=J_{0}(z_{1},z_{2})$ we have
\begin{equation}
\Phi (z_{1},z_{2})+\left| J_{2}(z_{1},z_{2})\right| \leq \Lambda
(z_{2})-\Lambda (z_{1})=:\Lambda (z_{2},z_{1})\text{ , for all }z_{2}\geq
z_{1}\text{.}  \label{sab2}
\end{equation}
by \eqref{eq:Jq(z)=-Hq(z) phi(z)}. Now define refining partitions
\begin{equation*}
-\infty =z_{0}(k)\leq z_{1}(k)\leq ....\leq z_{2^{k}}(k)=\infty \text{,}
\end{equation*}%
where
\begin{equation*}
z_{i}(k)=\inf \left\{ z:\Lambda (z)\geq \frac{\Lambda (\infty )i}{2^{k}}%
\right\} \text{ },\text{ }i=0,...,2^{k}-1\text{,}
\end{equation*}%
implying that%
\begin{equation*}
\Lambda (z_{i}(k))-\Lambda (z_{i-1}(k))\leq \frac{1.5}{2^{k}}\text{.}
\end{equation*}%
Define $i_{k}(z)$ by
\begin{equation*}
z_{i_{k}(z)}(k)\leq z\leq z_{i_{k}(z)+1}(k)\text{,}
\end{equation*}%
and note that the sequence $\left\{ z_{i_{k}(z)}(k)\right\} $ is increasing,
while $\left\{ z_{i_{k}(z)+1}(k)\right\} $ is decreasing; then write, for
some integer $K=K(l)$ to be determined later%
\begin{equation}  \label{eq:Sl(z) as sum}
S_{l}(z)=S_{l}(z_{i_{0}(z)}(0),z_{i_{1}(z)}(1))+S_{l}(z_{i_{1}(z)}(1),z_{i_{2}(z)}(2))+...+S_{l}(z_{i_{K}(z)}(K),z)%
\text{.}
\end{equation}%
The idea is to bound uniformly each term of \eqref{eq:Sl(z) as sum} by means
of (\ref{sab1}) and (\ref{sab2}). For the last term of \eqref{eq:Sl(z) as
sum} we have, by a simple algebraic manipulation
\begin{eqnarray*}
\left| S_{l}(z_{i_{K}(z)}(K),z)\right| &\leq &\sqrt{l}\left\{
\int_{S^{2}}\left\{ \mathbb{I}(z_{i_{K}}(K)\leq f_{l}(x)\leq
z_{i_{K}+1}(K))\right\} dx+4\pi \Phi (z_{i_{K}}(K),z_{i_{K}+1}(K))\right\} \\
&&+\sqrt{l}\Lambda (z_{i_{K}}(K),z_{i_{K}+1}(K))\left| h_{l;2}\right| \\
&\leq &\sqrt{l}\left| S_{l}(z_{i_{K}(z)}(K),z_{i_{K}(z)+1}(K))\right|
+2\times 4\pi \sqrt{l}\Lambda (z_{i_{K}}(K),z_{i_{K}+1}(K)) \\
&&+\sqrt{l}\Lambda (z_{i_{K}}(K),z_{i_{K}+1}(K))\left| h_{l;2}\right| \\
&\leq &\sqrt{l}\left| S_{l}(z_{i_{K}(z)}(K),z_{i_{K}(z)+1}(K))\right| +\sqrt{
l}\frac{8\pi \times 1.5}{2^{K}} \\
&&+\sqrt{l}\frac{1.5}{2^{K}}\left| h_{l;2}\right| \text{.}
\end{eqnarray*}

Using the latter together with \eqref{eq:Sl(z) as sum}, we may bound the
latter as
\begin{equation*}
\begin{split}
\Pr \left\{ \sup_{z}\left| S_{l}(z)\right| >\epsilon \right\} & \leq \Pr
\left( \left| S_{l}(z_{i_{0}(z)}(0),z_{i_{1}(z)}(1))\right| +...+\left|
S_{l}(z_{i_{K}(z)}(K),z_{i_{K}(z)+1}(K))\right| >\frac{\epsilon }{2}\right)
\\
& +\Pr \left( \sqrt{l}\frac{8\pi \times 1.5}{2^{K}}+\sqrt{l}\frac{1.5}{2^{K}}%
\left| h_{l;2}\right| >\frac{\epsilon }{2}\right) .
\end{split}%
\end{equation*}%
Since $\sum_{k=0}^{\infty }\varepsilon /(k+3)^{2}\leq \varepsilon /2,$ we
may bound
\begin{equation}
\Pr \left\{ \sup_{z}\left| S_{l}(z)\right| >\varepsilon \right\} \leq \Pr
\left\{ \sup_{z}\left| S_{l}(z_{i_{0}(z)}(0),z_{i_{1}(z)}(1))\right| >\frac{%
\varepsilon }{9}\right\}   \label{eq:Sl > eps sum}
\end{equation}%
\begin{equation*}
+\Pr \left\{ \sup_{z}\left| S_{l}(z_{i_{1}(z)}(1),z_{i_{2}(z)}(2))\right| >%
\frac{\varepsilon }{16}\right\} +....
\end{equation*}%
\begin{equation*}
+\Pr \left\{ \sup_{z}\left| S_{l}(z_{i_{K}(z)}(K),z_{i_{K}(z)+1}(K))\right| >%
\frac{\varepsilon }{(K+3)^{2}}\right\}
\end{equation*}%
\begin{equation*}
+\Pr \left\{ \sqrt{l}\frac{1.5}{2^{K}}\left| h_{l;2}\right| >\frac{%
\varepsilon }{2}-\frac{12\pi }{2^{K}}\sqrt{l}\right\} \text{.}
\end{equation*}%
Here, using the idea of refining partitions, we effectively reduced the
supremum over a continuous set ($\mathbb{R}$) to the supremum over the
possible nearest neighbors of $z$ in these partitions. In view of %
\eqref{sab1}, the Chebyshev inequality yields
\begin{eqnarray}
&&\Pr \left\{ \sup_{z}\left|
S_{l}(z_{i_{k}(z)}(k),z_{i_{k+1}(z)}(k+1))\right| >\frac{\varepsilon }{%
(k+3)^{2}}\right\}   \label{eq:Pr Sl Cheb} \\
&\leq &\sum_{i=0}^{2^{k+1}-1}\Pr \left\{ \left|
S_{l}(z_{i}(k+1),z_{i+1}(k+1))\right| >\frac{\varepsilon }{(k+3)^{2}}%
\right\}
\end{eqnarray}%
\begin{equation*}
\leq \frac{C\log l}{l}\frac{(k+3)^{4}}{\varepsilon ^{2}}%
\sum_{i=0}^{2^{k+1}-1}\left[ \Phi (z_{i+1}(k+1))-\Phi (z_{i}(k+1))\right]
\leq \frac{C\log l}{l}\frac{(k+3)^{4}}{\varepsilon ^{2}}\text{.}.
\end{equation*}%
We choose
\begin{equation}
K=\left\lfloor \log _{2}\left( \frac{48}{\varepsilon }\pi \sqrt{l}\right)
\right\rfloor +1\text{, }  \label{paci}
\end{equation}%
so that
\begin{equation*}
2^{K}\approx \frac{48\pi }{\varepsilon }\sqrt{l}\text{, }\sqrt{l}\frac{1.5}{%
2^{K}}\approx \frac{\varepsilon }{32\pi }\text{ , }\frac{\varepsilon }{2}-%
\frac{12\pi }{2^{K}}\sqrt{l}\approx \frac{\varepsilon }{4}\text{.}
\end{equation*}%
Hence we obtain
\begin{equation}
\begin{split}
& \Pr \left\{ \sqrt{l}\frac{1.5}{2^{K}}\left| h_{l;2}\right| >\frac{%
\varepsilon }{2}-\frac{12\pi }{2^{K}}\sqrt{l}\right\} \leq \Pr \left\{ \sqrt{%
l}\frac{1.5}{2^{K}}\left| h_{l;2}dx\right| >\frac{\varepsilon }{4}\right\}
\leq \Pr \left\{ \left| h_{l;2}dx\right| >8\pi \right\}  \\
& \leq \frac{1}{(8\pi )^{2}}l^{-1}E\left\{ l^{1/2}h_{l;2}\right\} ^{2}=\frac{%
1}{(8\pi )^{2}}l^{-1}Var\left\{ l^{1/2}\sum_{m}\left| a_{lm}\right|
^{2}\right\} \leq \frac{1}{(8\pi )^{2}}l^{-1}\text{.}
\end{split}
\label{eq:pr h > eps-}
\end{equation}%
To conclude the proof of \eqref{chain1}, we plug the bounds \eqref{eq:Pr Sl
Cheb} and \eqref{eq:pr h > eps-} into \eqref{eq:Sl > eps sum}; we obtain
\begin{eqnarray*}
\Pr \left\{ \sup_{z}\left| S_{l}(z)\right| >\varepsilon \right\}  &\leq &%
\frac{C\log l}{l}\sum_{k=0}^{K}\frac{(k+3)^{4}}{\varepsilon ^{2}}+\frac{%
l^{-1}}{(8\pi )^{2}} \\
&\leq &\frac{C^{\prime }\log l}{l}\frac{(K+3)^{5}}{5\varepsilon ^{2}}+\frac{%
l^{-1}}{(8\pi )^{2}} \\
&\leq &\frac{C^{\prime \prime }}{l}\left\{ \frac{\log ^{6}l}{\varepsilon ^{3}%
}+1\right\} \text{, some }C,C^{\prime },C^{\prime \prime }>0\text{.}
\end{eqnarray*}%
For the proof of (\ref{chain2}), we repeat exactly the same argument as
before, replacing $K$ in (\ref{paci}) by
\begin{equation*}
K^{\prime }=\left\{ \log _{2}\left( \frac{48\pi }{\varepsilon ^{1/3}}\sqrt{l}%
\right) \right\} \text{ , }2^{K^{\prime }}=\frac{48\pi }{\varepsilon ^{1/3}}%
\sqrt{l}\text{, }
\end{equation*}%
\begin{equation*}
\frac{\varepsilon }{2}-\frac{12\pi }{2^{K^{\prime }}}\sqrt{l}\geq \frac{%
\varepsilon }{2}-\frac{\varepsilon ^{1/3}}{4}\geq \frac{\varepsilon }{4}%
\text{ , because }\varepsilon \geq 1\text{.}
\end{equation*}%
It follows easily that%
\begin{eqnarray*}
&&\Pr \left\{ \sqrt{l}\frac{1.5}{2^{K}}\left| \int_{S^{2}}\left\{
f_{l}^{2}(x)-1\right\} dx\right| >\frac{\varepsilon }{2}-\frac{12\pi }{2^{K}}%
\sqrt{l}\right\}  \\
&\leq &\Pr \left\{ \frac{\varepsilon ^{1/3}}{32\pi }\left|
\int_{S^{2}}\left\{ f_{l}^{2}(x)-1\right\} dx\right| >\frac{\varepsilon }{4}%
\right\}  \\
&\leq &\Pr \left\{ \left| \int_{S^{2}}\left\{ f_{l}^{2}(x)-1\right\}
dx\right| >8\pi \varepsilon ^{2/3}\right\} \leq \frac{\varepsilon ^{-4/3}}{%
(8\pi )^{2}}l^{-1},
\end{eqnarray*}%
and%
\begin{eqnarray*}
\Pr \left\{ \sup_{z}\left| S_{l}(z)\right| >\varepsilon \right\}  &\leq &%
\frac{C\log l}{l}\frac{(K^{\prime }+3)^{5}}{5\varepsilon ^{2}}+C^{\prime }%
\frac{\varepsilon ^{-4/3}}{(8\pi )^{2}}l^{-1} \\
&\leq &\frac{C^{\prime \prime }}{l}\left\{ \frac{\log ^{6}l}{\varepsilon
^{7/3}}+\varepsilon ^{-4/3}\right\} \text{, some }C,C^{\prime },C^{\prime
\prime }>0\text{.}
\end{eqnarray*}
\end{proof}

\begin{remark}
The proof above is uses the integral
\begin{equation*}
\frac{1}{2}\int_{-\infty }^{z}\left| x^{2}-1\right| \phi (x)dx=
\end{equation*}%
\begin{equation*}
=\frac{1}{2}\int_{-\infty }^{z}\left( x^{2}-1\right) \phi (x)dx\mathbb{I}%
(z\leq -1)
\end{equation*}%
\begin{equation*}
+\left\{ \frac{1}{2}\int_{-\infty }^{-1}\left( x^{2}-1\right) \phi (x)dx%
\mathbb{+}\frac{1}{2}\int_{-1}^{z}\left( 1-x^{2}\right) \phi (x)dx\right\}
\mathbb{I}(-1<z\leq 1)
\end{equation*}%
\begin{equation*}
+\left\{ \frac{1}{2}\int_{-\infty }^{-1}\left( x^{2}-1\right) \phi (x)dx%
\mathbb{+}\frac{1}{2}\int_{-1}^{1}\left( 1-x^{2}\right) \phi (x)dx\frac{1}{2}%
\int_{1}^{z}\left( x^{2}-1\right) \phi (x)dx\right\} \mathbb{I}(z>1)
\end{equation*}%
\begin{equation*}
=\frac{1}{2\sqrt{2\pi }}\left\{ -z\exp (-\frac{z^{2}}{2})\mathbb{I}(z\leq
-1)+\left[ z\exp (-\frac{z^{2}}{2})+\exp (-\frac{1}{2})\right] \mathbb{I}%
(-1<z\leq 1)\right\}
\end{equation*}%
\begin{equation*}
+\frac{1}{2\sqrt{2\pi }}\left[ -z\exp (-\frac{z^{2}}{2})+4\exp (-\frac{1}{2})%
\right] \mathbb{I}(z>1),
\end{equation*}
which follows directly from
\begin{equation*}
\int_{-\infty }^{z}\left( x^{2}-1\right) \phi (x)dx\mathbb{I}(z\leq -1)=%
\left[ x\phi (x)\right] _{-\infty }^{z}=-\frac{z}{\sqrt{2\pi }}\exp (-\frac{%
z^{2}}{2})\mathbb{I}(z\leq -1),
\end{equation*}%
\begin{equation*}
\int_{-1}^{z}\left( 1-x^{2}\right) \phi (x)dx\mathbb{I}(z\leq 1)=\frac{1}{%
\sqrt{2\pi }}\left[ z\exp (-\frac{z^{2}}{2})+\exp (-\frac{1}{2})\right]
\mathbb{I}(z\leq 1),
\end{equation*}
and
\begin{equation*}
\int_{1}^{z}\left( x^{2}-1\right) \phi (x)dx\mathbb{I}(z>1)=\frac{1}{\sqrt{%
2\pi }}\left[ -z\exp (-\frac{z^{2}}{2})+\exp (-\frac{1}{2})\right] \mathbb{I}%
(z>1)\text{.}
\end{equation*}%
Also%
\begin{equation*}
\frac{1}{2}\int_{-\infty }^{\infty }\left| x^{2}-1\right| \phi (x)dx=\frac{%
4e^{-1/2}}{2\sqrt{2\pi }}\simeq 0.48394\text{.}
\end{equation*}
\end{remark}

\section{Proof of Theorem \ref{thm:Wl->Winf}}

In view of (\ref{important}), we can write%
\begin{eqnarray*}
W_{L}(z;r) &=&\frac{J_{2}(z)}{\sqrt{L}}\sum_{l=1}^{[Lr]}\left\{ \sqrt{l}%
\int_{S^{2}}\left\{ f_{l}^{2}(x)-1\right\} dx\right\} +\frac{1}{\sqrt{L}}%
\sum_{l=1}^{[Lr]}S_{l}(z) \\
&=&W_{AL}(z;r)+W_{BL}(z;r)\text{.}
\end{eqnarray*}%
We shall prove that, as $L\rightarrow \infty ,$
\begin{equation*}
W_{AL}(z;r)\Rightarrow W_{\infty }(z;r)\text{, }\sup_{z}%
\sup_{r}W_{BL}(z;r)=o_{p}(1)\text{.}
\end{equation*}

\subsection{Step 1: the proof that $W_{AL}(z;r)\Rightarrow W_{\infty }(z;r),$
as $L\rightarrow \infty $}

\begin{proof}
To prove convergence of the finite-dimensional distributions for $%
W_{AL}(z;r) $, it is enough to note that, for $r_{1}\leq r_{2}\leq r_{3}\leq
r_{4}$
\begin{eqnarray*}
EW_{AL}(z;r) &=&\frac{J_{2}(z)}{\sqrt{L}}\sum_{l=1}^{[Lr]}\sqrt{l}%
\int_{S^{2}}E\left\{ f_{l}^{2}(x)-1\right\} dx=0\text{,} \\
EW_{AL}(z_{1};r_{1})W_{AL}(z_{2};r_{2}) &=&\frac{J_{2}(z_{1})J_{2}(z_{2})}{L}%
\sum_{l=1}^{[Lr_{1}]}lEh_{l,2}^{2}\rightarrow J_{2}(z_{1})J_{2}(z_{2})r_{1}%
\text{,}
\end{eqnarray*}%
and
\begin{eqnarray*}
&&cum\left\{
W_{AL}(z_{1};r_{1}),W_{AL}(z_{2};r_{2}),W_{AL}(z_{3};r_{3}),W_{AL}(z_{4};r_{4}),\right\}
\\
&=&\frac{J_{2}(z_{1})J_{2}(z_{2})J_{3}(z_{1})J_{4}(z_{2})}{L^{2}}%
\sum_{l=1}^{[Lr_{1}]}l^{2}cum_{4}(h_{l,2})\rightarrow 0\text{, as }%
L\rightarrow \infty \text{.}
\end{eqnarray*}%
The multivariate extension is trivial. To establish tightness, in view of
the results from \cite{biwi} (see also \cite{mapi}), it is enough to prove
that, for all $r_{1}\leq r\leq r_{2},z_{1}\leq z\leq z_{2},$ $r=k/L,$ $%
k=1,...,L,$ some $C>0,$ there exist a finite measure $\mu :\mathcal{B}%
([0,1]\times \mathbb{R)\rightarrow R}^{+}$ such that%
\begin{eqnarray*}
&&E\left\{ \left[ W_{AL}(z_{2};r)-W_{AL}(z_{1};r_{1})\right] ^{2}\left[
W_{AL}(z_{2};r_{2})-W_{AL}(z_{1};r)\right] ^{2}\right\} \\
&=&E\left\{ \left[ \frac{J_{2}(z_{2})-J_{2}(z_{1})}{\sqrt{L}}%
\sum_{l=[Lr_{1}]}^{[Lr]}\sqrt{l}h_{l;2}\right] ^{2}\left[ \frac{%
J_{2}(z_{2})-J_{2}(z_{1})}{\sqrt{L}}\sum_{l=[Lr]+1}^{[Lr_{2}]}\sqrt{l}h_{l;2}%
\right] ^{2}\right\}
\end{eqnarray*}%
\begin{eqnarray}
&&  \notag \\
&\leq &C\left\{ \mu ([r_{1},r_{2}]\times \lbrack z_{1},z_{2}])\right\} ^{2},
\label{sat}
\end{eqnarray}%
\begin{eqnarray*}
&&E\left\{ \left[ W_{AL}(z;r_{2})-W_{AL}(z_{1};r_{1})\right] ^{2}\left[
W_{AL}(z_{2};r_{2})-W_{AL}(z;r_{1})\right] ^{2}\right\} \\
&=&E\left\{ \left[ \frac{J_{2}(z)-J_{2}(z_{1})}{\sqrt{L}}%
\sum_{l=[Lr_{1}]}^{[Lr_{2}]}\sqrt{l}h_{l;2}\right] ^{2}\left[ \frac{%
J_{2}(z_{2})-J_{2}(z)}{\sqrt{L}}\sum_{l=[Lr_{1}]}^{[Lr_{2}]}\sqrt{l}h_{l;2}%
\right] ^{2}\right\}
\end{eqnarray*}%
\begin{eqnarray}
&&  \notag \\
&\leq &C\left\{ \mu ([r_{1},r_{2}]\times \lbrack z_{1},z_{2}])\right\} ^{2}.
\label{sat2}
\end{eqnarray}%
In our case, both these bounds follow easily considering the measure%
\begin{equation*}
\mu ([r_{1},r_{2}]\times \lbrack z_{1},z_{2}]):=(r_{2}-r_{1})\left\{ \Lambda
_{2}(z_{2})-\Lambda _{1}(z_{1})\right\} \text{, }\mu (\mathbb{R\times R)<}%
2.25<\infty \text{.}
\end{equation*}%
Indeed, for (\ref{sat}), it is enough to exploit independence over $l$ to
show that%
\begin{eqnarray*}
&&E\left\{ \left[ \frac{J_{2}(z_{2})-J_{2}(z_{1})}{\sqrt{L}}%
\sum_{l=[Lr_{1}]}^{[Lr]}\sqrt{l}h_{l;2}\right] ^{2}\left[ \frac{%
J_{2}(z_{2})-J_{2}(z_{1})}{\sqrt{L}}\sum_{l=[Lr]+1}^{[Lr_{2}]}\sqrt{l}h_{l;2}%
\right] ^{2}\right\} \\
&\leq &\left[ \left\{ \frac{J_{2}(z_{2})-J_{2}(z_{1})}{\sqrt{L}}\right\}
^{2}\sum_{l=[Lr_{1}]}^{[Lr]}E\left\{ lh_{l;2}\right\} ^{2}\right] \left[
\left\{ \frac{J_{2}(z_{2})-J_{2}(z_{1})}{\sqrt{L}}\right\}
^{2}\sum_{l=[Lr]+1}^{[Lr_{2}]}E\left\{ lh_{l;2}^{2}\right\} \right] \\
&\leq &(r-r_{1})(r_{2}-r)\left\{ J_{2}(z_{2})-J_{2}(z_{1})\right\} ^{2}\leq
(r_{2}-r_{1})^{2}\left\{ \Lambda _{2}(z_{2})-\Lambda _{1}(z_{1})\right\}
^{2}.
\end{eqnarray*}%
For (\ref{sat2}), we have%
\begin{equation*}
E\left\{ \left[ \frac{J_{2}(z)-J_{2}(z_{1})}{\sqrt{L}}%
\sum_{l=[Lr_{1}]}^{[Lr_{2}]}\sqrt{l}h_{l;2}\right] ^{2}\left[ \frac{%
J_{2}(z_{2})-J_{2}(z)}{\sqrt{L}}\sum_{l=[Lr_{1}]}^{[Lr_{2}]}\sqrt{l}h_{l;2}%
\right] ^{2}\right\}
\end{equation*}%
\begin{eqnarray*}
&=&\frac{\left\{ J_{2}(z)-J_{2}(z_{1})\right\} ^{2}\left\{
J_{2}(z)-J_{2}(z_{1})\right\} ^{2}}{L^{2}}\sum_{l=[Lr_{1}]}^{[Lr_{2}]}E\left%
\{ \sqrt{l}h_{l;2}\right\} ^{4} \\
&&+6\frac{\left\{ J_{2}(z)-J_{2}(z_{1})\right\} ^{2}\left\{
J_{2}(z)-J_{2}(z_{1})\right\} ^{2}}{L^{2}}%
\sum_{l_{1}<l_{2}=[Lr_{1}]}^{[Lr_{2}]}E\left\{ \sqrt{l_{1}}%
h_{l_{1};2}\right\} ^{2}E\left\{ \sqrt{l_{2}}h_{l_{2};2}\right\} ^{2}.
\end{eqnarray*}%
Now $E\left\{ \sqrt{l}h_{l;2}\right\} ^{2},E\left\{ \sqrt{l}h_{l;2}\right\}
^{4}\leq C$ uniformly w.r.t. $l;$ note also that $L^{-1}\leq (r_{2}-r_{1}),$
for all $r_{1},r_{2}$ in $\left\{ 1/L,2/L,...,1\right\} ,$ whence%
\begin{eqnarray*}
&&\frac{\left\{ J_{2}(z)-J_{2}(z_{1})\right\} ^{2}\left\{
J_{2}(z)-J_{2}(z_{1})\right\} ^{2}}{L^{2}}\sum_{l=[Lr_{1}]}^{[Lr_{2}]}E\left%
\{ \sqrt{l}h_{l;2}\right\} ^{4} \\
&\leq &C(r_{2}-r_{1})^{2}\left\{ \Lambda _{2}(z_{2})-\Lambda
_{1}(z_{1})\right\} ^{2}, \\
&&\frac{\left\{ J_{2}(z)-J_{2}(z_{1})\right\} ^{2}\left\{
J_{2}(z)-J_{2}(z_{1})\right\} ^{2}}{L^{2}}%
\sum_{l_{1}<l_{2}=[Lr_{1}]}^{[Lr_{2}]}E\left\{ \sqrt{l_{1}}%
h_{l_{1};2}\right\} ^{2}E\left\{ \sqrt{l_{2}}h_{l_{2};2}\right\} ^{2} \\
&\leq &(r_{2}-r_{1})^{2}\left\{ \Lambda _{2}(z_{2})-\Lambda
_{1}(z_{1})\right\} ^{2},
\end{eqnarray*}%
as needed.
\end{proof}

\subsection{Step 2: the proof that $\sup_{z}\sup_{r}W_{BL}(z;r)=o_{p}(1)$ ,
as $L\rightarrow \infty $ .}

\begin{proof}
For fixed $z,r$ it is enough to show that%
\begin{equation*}
\sup_{0\leq r\leq 1}\sup_{-\infty <z<\infty }\left| W_{BL}(z;r)\right|
=o_{p}(1)\text{ , as }L\rightarrow \infty \text{.}
\end{equation*}%
To establish this result, we note that%
\begin{eqnarray*}
E\left\{ \sup_{0\leq r\leq 1}\sup_{-\infty <z<\infty }\left|
W_{BL}(z;r)\right| \right\} &\leq &\frac{1}{\sqrt{L}}E\left\{ \sup_{0\leq
r\leq 1}\sum_{l=1}^{[Lr]}\sup_{-\infty <z<\infty }\left| S_{l}(z)\right|
\right\} \\
&\leq &\frac{1}{\sqrt{L}}\sum_{l=1}^{L}E\left\{ \sup_{-\infty <z<\infty
}\left| S_{l}(z)\right| \right\} \\
&\leq &\frac{C}{\sqrt{L}}\sum_{l=1}^{L}\frac{\log ^{6}l}{l}=C\frac{\log ^{7}L%
}{\sqrt{L}}\rightarrow 0\text{ , as }L\rightarrow \infty \text{ ,}
\end{eqnarray*}%
because
\begin{eqnarray*}
E\left\{ \sup_{-\infty <z<\infty }\left| S_{l}(z)\right| \right\} &\leq &1+%
\frac{C\log l}{l}\int_{1}^{\infty }\left\{ \varepsilon ^{-4/3}+\frac{\log
^{5}l}{\varepsilon ^{3}}\right\} d\varepsilon \\
&\leq &C^{\prime }\frac{\log ^{6}l}{l}\text{ .}
\end{eqnarray*}%
\begin{eqnarray*}
W_{BL}(z;r) &=&O_{p}\left( \left[ EW_{BL}^{2}(z;r)\right] ^{1/2}\right)
=O_{p}\left( \frac{1}{L}\sum_{l=1}^{[Lr]}ES_{l}^{2}(z)\right) \\
&=&O_{p}\left( \frac{1}{L}\sum_{l=1}^{[Lr]}\frac{\log l}{l}\right) =o_{p}(1).
\end{eqnarray*}
\end{proof}

\appendix

\section{Background on Wigner Coefficients}

Throughout this paper, we made a heavy use of Wigner's $3j$ coefficients. In
this appendix, we review briefly some of their features and provide some
results on their asymptotic properties. We refer to \cite{VilKly}, \cite{VMK}
and \cite{BieLou} for a much more detailed discussion, in particular
concerning the relationships with the quantum theory of angular momentum and
group representation properties of $SO(3).$

We start from the analytic expression (valid for $m_{1}+m_{2}+m_{3}=0,$ see %
\cite{VMK}, expression (8.2.1.5))
\begin{align*}
\left(
\begin{array}{ccc}
l_{1} & l_{2} & l_{3} \\
m_{1} & m_{2} & m_{3}%
\end{array}%
\right) & :=(-1)^{l_{1}+m_{1}}\sqrt{2l_{3}+1}\left[ \frac{%
(l_{1}+l_{2}-l_{3})!(l_{1}-l_{2}+l_{3})!(l_{1}-l_{2}+l_{3})!}{%
(l_{1}+l_{2}+l_{3}+1)!}\right] ^{1/2} \\
& \times \left[ \frac{(l_{3}+m_{3})!(l_{3}-m_{3})!}{%
(l_{1}+m_{1})!(l_{1}-m_{1})!(l_{2}+m_{2})!(l_{2}-m_{2})!}\right] ^{1/2} \\
& \times \sum_{z}\frac{(-1)^{z}(l_{2}+l_{3}+m_{1}-z)!(l_{1}-m_{1}+z)!}{%
z!(l_{2}+l_{3}-l_{1}-z)!(l_{3}+m_{3}-z)!(l_{1}-l_{2}-m_{3}+z)!}\text{,}
\end{align*}%
where the summation runs over all $z$'s such that the factorials are
non-negative. This expression becomes much neater for $m_{1}=m_{2}=m_{3}=0,$
where we have%
\begin{equation*}
\left(
\begin{array}{ccc}
l_{1} & l_{2} & l_{3} \\
0 & 0 & 0%
\end{array}%
\right) =
\end{equation*}%
\begin{equation}
\left\{
\begin{array}{c}
0\text{ , for }l_{1}+l_{2}+l_{3}\text{ odd} \\
(-1)^{\frac{l_{1}+l_{2}-l_{3}}{2}}\frac{\left[ (l_{1}+l_{2}+l_{3})/2\right] !%
}{\left[ (l_{1}+l_{2}-l_{3})/2\right] !\left[ (l_{1}-l_{2}+l_{3})/2\right] !%
\left[ (-l_{1}+l_{2}+l_{3})/2\right] !}\left\{ \frac{%
(l_{1}+l_{2}-l_{3})!(l_{1}-l_{2}+l_{3})!(-l_{1}+l_{2}+l_{3})!}{%
(l_{1}+l_{2}+l_{3}+1)!}\right\} ^{1/2}\text{ } \\
\text{for }l_{1}+l_{2}+l_{3}\text{ even}%
\end{array}%
\right. .  \label{appe}
\end{equation}%
Some of the properties to follow become neater when expressed in terms of
the so-called Clebsch-Gordan coefficients, which are defined by by the
identities (see \cite{VMK}, Chapter 8)%
\begin{eqnarray}
\left(
\begin{array}{ccc}
l_{1} & l_{2} & l_{3} \\
m_{1} & m_{2} & -m_{3}%
\end{array}%
\right)  &=&(-1)^{l_{3}+m_{3}}\frac{1}{\sqrt{2l_{3}+1}}%
C_{l_{1}-m_{1}l_{2}-m_{2}}^{l_{3}m_{3}}  \label{clewig1} \\
C_{l_{1}m_{1}l_{2}m_{2}}^{l_{3}m_{3}} &=&(-1)^{l_{1}-l_{2}+m_{3}}\sqrt{%
2l_{3}+1}\left(
\begin{array}{ccc}
l_{1} & l_{2} & l_{3} \\
m_{1} & m_{2} & -m_{3}%
\end{array}%
\right) \text{.}  \label{clewig2}
\end{eqnarray}%
We have the following orthonormality conditions:
\begin{eqnarray}
\sum_{m_{1},m_{2}}C_{l_{1}m_{1}l_{2}m_{2}}^{lm}C_{l_{1}m_{1}l_{2}m_{2}}^{l^{%
\prime }m^{\prime }} &=&\delta _{l}^{l\prime }\delta _{m}^{m\prime },
\label{ortho1} \\
\sum_{l,m}C_{l_{1}m_{1}l_{2}m_{2}}^{lm}C_{l_{1}m_{1}^{\prime
}l_{2}m_{2}^{\prime }}^{lm} &=&\delta _{m_{1}}^{m_{1}^{\prime }}\delta
_{m_{2}}^{m_{2}^{\prime }}\text{.}  \label{ortho2}
\end{eqnarray}%
Now recall the general formula (\cite{VMK}, eqs. (5.6.2.12-13))%
\begin{equation*}
\int_{S^{2}}Y_{l_{1}m_{1}}(x)...Y_{l_{n}m_{n}}(x)dx
\end{equation*}%
\begin{eqnarray*}
&=&\sqrt{\frac{4\pi }{2l_{n}+1}}\sum_{L_{1}...L_{n-3}}\sum_{M_{1}...M_{n-3}}%
\left[
C_{l_{1}m_{1}l_{2}m_{2}}^{L_{1}M_{1}}C_{L_{1}M_{1}l_{3}m_{3}}^{L_{2}M_{2}}...C_{L_{n-3}M_{n-3}l_{n-1}m_{n-1}}^{l_{n}m_{n}}\right.
\\
&&\times \left. \sqrt{\frac{\prod_{i=1}^{n-1}(2l_{i}+1)}{(4\pi )^{n-1}}}%
\left\{
C_{l_{1}0l_{2}0}^{L_{1}0}C_{L_{1}0l_{3}0}^{L_{2}0}...C_{L_{n-3}0l_{n-1}0}^{l_{n}0}\right\} %
\right] \text{.}
\end{eqnarray*}%
Hence we have%
\begin{eqnarray*}
\int_{0}^{1}P_{l}^{n}(t)dt &=&\sqrt{\frac{(4\pi )^{n-2}}{(2l+1)^{n}}}%
\int_{S^{2}}Y_{l0}^{n}(x)dx=\sum_{L_{1}...L_{n-3}}\left\{
C_{l_{1}0l_{2}0}^{L_{1}0}C_{L_{1}0l_{3}0}^{L_{2}0}...C_{L_{n-3}0l_{n-1}0}^{l_{n}0}\right\} ^{2}
\\
&=&\frac{1}{2l+1}\sum_{L_{1}...L_{n-3}}\left\{
C_{l_{1}0l_{2}0...l_{n-1}0}^{L_{1}L_{2}...l0}\right\} ^{2},
\end{eqnarray*}%
in the notation of \cite{MaPeBer}, \cite{MaPeJMVA}. Special cases are%
\begin{eqnarray*}
\int_{0}^{1}P_{l}^{3}(t)dt &=&\frac{1}{2}\int_{0}^{\pi }P_{l}^{3}(\cos
\vartheta )d\cos \vartheta  \\
&=&\sqrt{\frac{4\pi }{(2l+1)^{3}}}\int_{0}^{2\pi }\int_{0}^{\pi
}Y_{l0}^{3}(\vartheta ,\varphi )\sin \vartheta d\vartheta d\varphi  \\
&=&\frac{1}{2l+1}\left\{ C_{l0l0}^{l0}\right\} ^{2}=\left(
\begin{array}{ccc}
l & l & l \\
0 & 0 & 0%
\end{array}%
\right) ^{2},
\end{eqnarray*}%
\begin{eqnarray}
\int_{0}^{1}P_{l}^{4}(t)dt &=&\frac{1}{2l+1}\sum_{L=0}^{2l}\left\{
C_{l0l0l0}^{Ll0}\right\} ^{2}=\frac{1}{2l+1}\sum_{L=0}^{2l}\left\{
C_{l0l0}^{L0}C_{L0l0}^{l0}\right\} ^{2}  \notag \\
&=&\frac{1}{2l+1}\sum_{L=0}^{2l}(2L+1)(2l+1)\left(
\begin{array}{ccc}
l & l & L \\
0 & 0 & 0%
\end{array}%
\right) ^{4}  \notag \\
&=&\sum_{L=0}^{2l}(2L+1)\left(
\begin{array}{ccc}
l & l & L \\
0 & 0 & 0%
\end{array}%
\right) ^{4},  \label{appe2}
\end{eqnarray}%
compare \cite{VMK}, equation (8.9.4.20). The following result is certainly
known, but we failed to locate a reference and we provide a proof for
completeness.

\begin{lemma}
\label{cgbou} 1) As $l\rightarrow \infty $%
\begin{equation*}
\lim_{l\rightarrow \infty }l^{2}\left(
\begin{array}{ccc}
l & l & l \\
0 & 0 & 0%
\end{array}%
\right) ^{2}=\frac{2}{\pi \sqrt{3}}\simeq 0.367\text{.}
\end{equation*}%
2) For all $l=1,2,..,$ \ and even $L=2,...,2l-2$ we have%
\begin{equation*}
\left(
\begin{array}{ccc}
l & l & L \\
0 & 0 & 0%
\end{array}%
\right) ^{2}=\gamma _{lL}\times \frac{2}{\pi }\times \frac{1}{%
L(2l-L)^{1/2}(2l+L)^{1/2}}\text{,}
\end{equation*}%
where%
\begin{equation*}
0.596=1.09^{-6}\leq \gamma _{lL}\leq 1.09^{5}=1.539\text{.}
\end{equation*}%
3) For $L=0,2l$ we have
\begin{equation*}
\left(
\begin{array}{ccc}
l & l & 0 \\
0 & 0 & 0%
\end{array}%
\right) ^{2}=\frac{1}{2l+1}\text{ , }\left(
\begin{array}{ccc}
l & l & 2l \\
0 & 0 & 0%
\end{array}%
\right) ^{2}=\frac{\sqrt{2}}{\sqrt{\pi }(4l+1)\sqrt{l}}\left\{
1+O(l^{-1})\right\} \text{.}
\end{equation*}
\end{lemma}

\begin{proof}
1) Note that (see (\ref{appe}) and \cite{VMK}, equation (8.5.2.32))%
\begin{equation*}
\left(
\begin{array}{ccc}
l & l & l \\
0 & 0 & 0%
\end{array}%
\right) =\frac{(-1)^{l/2}(3l/2)!}{\left[ (l/2)!\right] ^{3}}\left[ \frac{%
\left[ l!\right] ^{3}}{(3l+1)!}\right] ^{1/2}.
\end{equation*}%
Hence, recalling Stirling's formula%
\begin{equation*}
l!=\sqrt{2\pi }(l)^{l+1/2}\exp (-l)+O(l^{-1})\text{,}
\end{equation*}%
we obtain%
\begin{eqnarray*}
\lim_{l\rightarrow \infty }l^{2}\left(
\begin{array}{ccc}
l & l & l \\
0 & 0 & 0%
\end{array}%
\right) ^{2} &=&\frac{\sqrt{(2\pi )^{5}}}{\sqrt{(2\pi )^{7}}}%
\lim_{l\rightarrow \infty }l^{2}\frac{(3l/2)^{3l+1}e^{-3l}}{%
(l/2)^{3l+3}e^{-3l}}\left[ \frac{l^{3l+3/2}e^{-3l}}{(3l+1)^{3l+3/2}e^{-3l-1}}%
\right] \\
&=&\frac{e}{2\pi }\lim_{l\rightarrow \infty }l^{2}\frac{(3l/2)^{3l+1}}{%
(l/2)^{3l+3}}\left[ \frac{l^{3l+3/2}}{(3l+1)^{3l+3/2}}\right] \\
&=&\frac{2e}{\pi }\lim_{l\rightarrow \infty }3^{3l+1}\left[ \frac{l^{3l+3/2}%
}{(3l+1)^{3l+3/2}}\right] \\
&=&\frac{6e}{3^{3/2}\pi }\lim_{l\rightarrow \infty }3^{3l}\left[ \frac{l^{3l}%
}{(3l+1)^{3l}}\right] =\frac{6e}{3^{3/2}\pi }\lim_{l\rightarrow \infty }%
\left[ \frac{1}{(1+l^{-1})^{3l}}\right] \\
&=&\frac{2}{\pi \sqrt{3}}\simeq 0.367\text{.}
\end{eqnarray*}%
2) From \cite{VMK}, equation (8.5.2.32) we have%
\begin{equation*}
\left(
\begin{array}{ccc}
l & l & L \\
0 & 0 & 0%
\end{array}%
\right) ^{2}=\left\{ \frac{(l+L/2)!}{((L/2)!)^{2}(l-L/2)!}\right\} ^{2}\frac{%
((L)!)^{2}((2l-L)!)}{(2l+L+1)!}\text{.}
\end{equation*}%
We use repeatedly Stirling's formula, for $n=1,2,...$
\begin{equation*}
1<\exp (\frac{1}{12n+1})\leq \frac{n!}{n^{n}e^{-n}\sqrt{2\pi n}}\leq \exp (%
\frac{1}{12n})\leq 1.09\text{;}
\end{equation*}%
also we write $a_{n}\asymp b_{n}$ for sequences such that $%
a_{n}/b_{n},b_{n}/a_{n}=O(1)$ . Hence we have%
\begin{equation*}
\left(
\begin{array}{ccc}
l & l & L \\
0 & 0 & 0%
\end{array}%
\right) ^{2}=\left\{ \frac{(l+L/2)!}{((L/2)!)^{2}(l-L/2)!}\right\} ^{2}\left[
\frac{(L!)^{2}(2l-L)!}{(2l+L+1)!}\right]
\end{equation*}%
\begin{equation*}
=\frac{\gamma _{lL}}{2\pi }\frac{e^{-2l-L}}{e^{-2L}e^{-2l+L}}\frac{%
e^{-2L}e^{-2l+L}}{e^{-2l-L-1}}\frac{(l+L/2)^{2l+2L+1}}{%
(L/2)^{2L+2}(l-L/2)^{2l-2L+1}}\frac{(L)^{2L+1}(2l-L)^{2l-L+1/2}}{%
(2l+L+1)^{2l+L+3/2}}
\end{equation*}%
\begin{equation*}
=\frac{\gamma _{lL}\times e}{2\pi }\frac{(l+L/2)^{2l+2L+1}}{%
(L/2)^{2L+2}(l-L/2)^{2l-2L+1}}\frac{(L)^{2L+1}(2l-L)^{2l-L+1/2}}{%
(2l+L+1)^{2l+L+3/2}}\text{,}
\end{equation*}%
where%
\begin{equation*}
0.596=1.09^{-6}<\exp (-\frac{6}{13})\leq \gamma _{lL}\leq \exp (\frac{5}{12}%
)<1.09^{5}=1.539\text{.}
\end{equation*}%
Now%
\begin{eqnarray*}
&&\frac{(l+L/2)^{2l+2L+1}}{(L/2)^{2L+2}(l-L/2)^{2l-2L+1}}\frac{%
(L)^{2L+1}(2l-L)^{2l-L+1/2}}{(2l+L+1)^{2l+L+3/2}} \\
&=&\frac{2^{2L+2}2^{2l-L+1/2}}{2^{2l+L+3/2}}\frac{(l+L/2)^{2l+L+1}}{%
L^{2L+2}(l-L/2)^{2l-L+1}}\frac{L^{2L+1}(l-L/2)^{2l-L+1/2}}{%
(l+L/2+1/2)^{2l+L+3/2}} \\
&=&\frac{2}{L(l-L/2)^{1/2}(l+L/2)^{1/2}}\frac{1}{(1+\frac{1}{2(l+L/2)}%
)^{2l+L+3/2}}=\frac{2e^{-1}}{L(l-L/2)^{1/2}(l+L/2)^{1/2}}\text{,}
\end{eqnarray*}%
whence the proof of the first is completed. \newline
3) The first part is equation (8.5.1.1) from (\cite{VMK}). For the second
part, it is sufficient to recall from \cite{VMK} equation (8.5.2.33) to
deduce that%
\begin{eqnarray*}
\left(
\begin{array}{ccc}
l & l & 2l \\
0 & 0 & 0%
\end{array}%
\right) ^{2} &=&\frac{1}{4l+1}\left\{ \frac{(2l)!}{l!l!}\left[ \frac{%
(2l)!(2l)!}{(4l)!}\right] ^{1/2}\right\} ^{2} \\
&=&\frac{1}{4l+1}\left\{ \frac{(2l)^{2l+1/2}}{\sqrt{2\pi }l^{l+1/2}l^{l+1/2}}%
\left[ \frac{\sqrt{2\pi }(2l)^{2l+1/2}(2l)^{2l+1/2}}{(4l)^{4l+1/2}}\right]
^{1/2}\right\} ^{2} \\
&=&\frac{\sqrt{2}}{\sqrt{\pi }(4l+1)\sqrt{l}}\text{,}
\end{eqnarray*}%
as claimed.
\end{proof}

\begin{remark}
Lemma \ref{cgbou} implies immediately
\begin{eqnarray*}
\left(
\begin{array}{ccc}
l & l & L \\
0 & 0 & 0%
\end{array}%
\right) ^{4} &\leq &\frac{1}{L^{2}\left\{ 4l^{2}-L^{2}\right\} }\text{ , for
}L=2,4....,2l-2 \\
&=&O(l^{-3})\text{ , for }L=2l\text{.}
\end{eqnarray*}%
Also, for $L=l$ we obtain%
\begin{equation*}
\frac{1}{L(l-L/2)^{1/2}(l+L/2)^{1/2}}=\frac{1}{l(l/2)^{1/2}(3l/2)^{1/2}}=%
\frac{2}{\sqrt{3}l^{2}}\text{,}
\end{equation*}%
leading to the special case%
\begin{equation*}
\lim_{l\rightarrow \infty }\left[ \frac{2}{\pi \sqrt{3}l^{2}}\right]
^{-1}\left(
\begin{array}{ccc}
l & l & l \\
0 & 0 & 0%
\end{array}%
\right) ^{2}=1\text{,}
\end{equation*}%
see also \cite{VMK}, equation (8.9.4.20) for different asymptotic
approximations.
\end{remark}

\section{Proof of Lemma \ref{quattro}}

\label{apx:proof lem Leg mom}

\begin{proof}
First, we note that, since for every $t\in \lbrack -1,1]$, $|P_{l}(t)|\leq 1$%
, it is sufficient to prove the statement of Lemma \ref{quattro} for $q=4$.
In this case, we have
\begin{equation*}
\begin{split}
E\left[ h_{l;4}^{2}\right] & =4!\int_{0}^{1}P_{l}^{4}(t)dt=12\int_{0}^{\pi
}P_{l}^{4}(\cos \vartheta )d\cos \vartheta \\
& =\frac{12}{4\pi }\int_{0}^{2\pi }\int_{0}^{\pi }P_{l}^{4}(\cos \vartheta
)d\cos \vartheta d\varphi =\frac{12}{4\pi }\sqrt{\left\{ \frac{4\pi }{2l+1}%
\right\} ^{4}}\int_{S^{2}}Y_{l0}^{4}(x)dx.
\end{split}%
\end{equation*}%
Now from \eqref{appe2} and Lemma \ref{cgbou}
\begin{equation*}
\begin{split}
\int_{0}^{2\pi }& \int_{0}^{\pi }P_{l}^{4}(\cos \vartheta )d\cos \vartheta
d\varphi =4\pi \int_{0}^{1}P_{l}^{4}(t)dt=\sum_{L=0}^{2l}(2L+1)\left(
\begin{array}{ccc}
l & l & L \\
0 & 0 & 0%
\end{array}%
\right) ^{4} \\
& \\
& \leq \frac{1}{(2l+1)^{2}}+\sum_{L=2}^{2L-2}(2L+1)\frac{1}{L^{2}\left\{
4l^{2}-L^{2}\right\} }+O(l^{-3}).
\end{split}%
\end{equation*}%
For $L\leq l$
\begin{equation*}
\sum_{L=2}^{l}(2L+1)\frac{1}{L^{2}\left\{ 4l^{2}-L^{2}\right\} }\leq \frac{C%
}{l^{2}}\sum_{L=2}^{l}\frac{1}{L}\leq C\frac{\log l}{l^{2}},
\end{equation*}%
while for $l\leq L<2l$%
\begin{equation*}
\sum_{L=l}^{2l-2}(2L+1)\frac{1}{L^{2}(2l+L)\left\{ 2l-L\right\} }\leq \frac{C%
}{2l^{2}}\sum_{L=l}^{2l-2}\frac{1}{\left\{ 2l-L\right\} }\leq C\frac{\log l}{%
l^{2}}\text{.}
\end{equation*}%
Note also that%
\begin{equation*}
\sum_{L=2}^{2l-2}(2L+1)\frac{1}{L^{2}\left\{ 4l^{2}-L^{2}\right\} }\geq
\frac{1}{4l^{2}}\sum_{L=2}^{2l-2}(2L+1)\frac{1}{L^{2}}\geq C^{\prime }\frac{%
\log l}{l^{2}},
\end{equation*}%
so this order cannot be improved. Also, using the second part of Lemma \ref%
{cgbou} we have easily
\begin{equation*}
(4l+1)\left(
\begin{array}{ccc}
l & l & 2l \\
0 & 0 & 0%
\end{array}%
\right) ^{4}=\frac{2}{\pi (4l+1)l}\left\{ 1+O(l^{-1})\right\} =O\left(
l^{-2}\right) .
\end{equation*}
\end{proof}

\end{document}